\RequirePackage{fix-cm}
\documentclass[11pt,letterpaper]{amsart}

\pdfoutput=1

\usepackage[obeyclassoptions]{standalone}

\usepackage[margin=1in]{geometry}
\usepackage[dvipsnames]{xcolor}

\usepackage{dsfont}
\usepackage{kpfonts}
\usepackage[T1]{fontenc}
\usepackage[utf8]{inputenc}

\usepackage{microtype}

\usepackage{mathtools} 
\usepackage[foot]{amsaddr}

\usepackage{amssymb}
\usepackage{mathtools} 

\usepackage{booktabs}

\usepackage{enumitem}
\usepackage{float,subcaption}

\usepackage{tikz,tkz-graph}

\usepackage[colorlinks=true, pdfstartview=FitV, linkcolor=blue, citecolor=teal, urlcolor=magenta]{hyperref}
\usepackage{bookmark}

\newtheorem{proposition}{Proposition}
\newtheorem{theorem}[proposition]{Theorem}
\newtheorem{lemma}[proposition]{Lemma}
\newtheorem{corollary}[proposition]{Corollary}

\newtheorem{conjecture}[proposition]{Conjecture}
\theoremstyle{remark}

\theoremstyle{definition}
\newtheorem{definition}[proposition]{Definition}

\numberwithin{equation}{section}
\numberwithin{proposition}{section}
\numberwithin{figure}{section}
\numberwithin{table}{section}

\newcommand{\E}{\mathbb{E}}
\renewcommand{\Pr}{\mathbb{P}}
\newcommand{\Var}{\mathbb{V}}
\renewcommand{\epsilon}{\varepsilon}

\newcommand*{\cI}{\mathcal{I}}
\newcommand{\geCOUNT}{\ge_{\mathrm{COUNT}}}
\newcommand{\gePART}{\ge_{\mathrm{PART}}}
\newcommand{\geCOEF}{\ge_{\mathrm{COEF}}}
\newcommand{\geOCC}{\ge_{\mathrm{OCC}}}
\newcommand{\geMAX}{\ge_{\mathrm{MAX}}}
\newcommand{\geFV}{\ge_{\mathrm{FV}}}
\newcommand{\geVAR}{\ge_{\mathrm{VAR}}}

\begin{document}

\author{Ewan Davies}
\address{Department of Computer Science, Colorado State University, USA}
\email{\href{mailto:ewan.davies@colostate.edu}{ewan.davies@colostate.edu}}

\author{Juspreet Singh Sandhu}
\address{Department of Computer Science, UC Santa Cruz, USA}
\email{\href{mailto:jsinghsa@ucsc.edu}{jsinghsa@ucsc.edu}}

\author{Brian Tan}
\address{Department of Computer Science, Colorado State University, USA.} \email{\href{mailto:brian.tan@colostate.edu}{brian.tan@colostate.edu}}


\title{On expectations and variances in the hard-core model on bounded degree graphs}

\begin{abstract}
    We extend the study of the occupancy fraction of the hard-core model in two novel directions. One direction gives a tight lower bound in terms of individual vertex degrees, extending work of Sah, Sawhney, Stoner and Zhao which bounds the partition function. The other bounds the variance of the size of an independent set drawn from the model, which is strictly stronger than bounding the occupancy fraction. 

    In the setting of triangle-free graphs, we make progress on a recent conjecture of Buys, van den Heuvel and Kang on extensions of Shearer's classic bounds on the independence number to the occupancy fraction of the hard-core model.

    Sufficiently strong lower bounds on both the expectation and the variance in triangle-free graphs have the potential to improve the known bounds on the off-diagonal Ramsey number $R(3,t)$, and to shed light on the algorithmic barrier one observes for independent sets in sparse random graphs.
\end{abstract}

\maketitle

\thispagestyle{empty}

\pagenumbering{arabic}

\section{Introduction}

We consider a distribution over independent sets in graphs known as the \emph{hard-core model} in which a set $I$ is chosen with probability proportional to $\lambda^{|I|}$ for some positive real parameter $\lambda$. 
This is a remarkably versatile distribution that generalizes the uniform one (which corresponds to $\lambda=1$), and that enjoys a number of important properties relating to conditional independence and entropy. 

In extremal graph theory, the hard-core model can be used to understand the minimum and maximum number of independent sets and matchings in regular graphs~\cite{DJPR17b}, the average size of independent sets in triangle-free graphs~\cite{DJPR17a}, and provides bounds on certain off-diagonal Ramsey numbers and related quantities~\cite{She95,Alo96,DKPS20a}. 
Some of these applications give, asymptotically, the best known results.
The conditional independence properties of the hard-core model can be exploited in probabilistic and algorithmic methods for graph coloring~\cite{Mol19,Ber19,BKNP22,DKPS20a,DKPS20b,DJKP20,DdKP21,Dha25}.

We study various extremal problems related to the number of independent sets in a graph and the average size of an independent set from the hard-core model. Extremal bounds in this context have a long history and a few surprising applications, e.g.~\cite{C87,A91,She95,Alo96,CDFKPY22}. 
Such extremal results determine a computational threshold for the problem of approximating the number of independent sets of a given size, and analogous thresholds in other spin systems~\cite{DP23,CDKP22,DL24}.
We also begin the systematic study of extremal behavior of the variance of the size of an independent set from the model which has, to our knowledge, escaped significant attention thus far.

Our main motivation is the so-called algorithmic barrier that one observes for independent sets in random $d$-regular graphs and binomial random graphs of constant average degree $d$. 
Here, for fixed $d$ it is well-known that independent sets of density asymptotically $2\log(d)/d$ exist (see e.g.~\cite{DJPR17a,BvK25} and the citations therein).
But the best-known efficient algorithms only find independent sets of density asymptotically half of this, which is tight for ``local'' algorithms~\cite{RV17}. 
Analysis of the hard-core model~\cite{BvK25,LMR+24a} can be used to demonstrate that certain algorithmic approaches hit the apparent barrier at density $\log(d)/d$, and our work is motivated by understanding and potentially overcoming the barrier.

Let $G$ be a graph on $n$ vertices, and let $\cI(G)$ be the set of independent sets in $G$, including the empty set. The \emph{partition function} of the hard-core model on $G$ is 
\[ Z_G(\lambda) = \sum_{I\in \cI(G)}\lambda^{|I|}, \]
and the \emph{hard-core model} on $G$ is the measure $\mu_{G,\lambda} : \cI(G)\to [0,1]$ given by 
\[ \mu_{G,\lambda}(I) = \frac{\lambda^{|I|}}{Z_G(\lambda)}. \]
Borrowing some physics terminology and inventing some of our own, we define the \emph{free energy density}
\[ F_G(\lambda) = \frac{1}{n}\log Z_G(\lambda), \]
the \emph{occupancy fraction}
\[ E_G(\lambda) = \lambda \frac{\partial}{\partial \lambda} F_G(\lambda), \]
and the \emph{variance fraction}
\[ V_G(\lambda) = \lambda \frac{\partial}{\partial \lambda} E_G(\lambda). \]
It is not too hard to see that 
\[ n E_G(\lambda) = \E_{I\sim \mu_{G,\lambda}}|I| \] 
is the expected number of vertices of $G$ occupied by a random independent set from the hard-core model, and similarly that 
\[ n V_G(\lambda) = \Var_{I\sim \mu_{G,\lambda}} |I| = \E_{I\sim \mu_{G,\lambda}}[|I|^2] - \E_{I\sim \mu_{G,\lambda}}[|I|]^2 \] 
is the variance of the size of an independent set from the model. 
The choice to normalize by $1/n$ is natural for the free energy $\log Z_G(\lambda)$ and expectation $\E_{I\sim \mu_{G,\lambda}}|I|$.
In an $n$-vertex graph the size of an independent set is supported on the interval $[0,n]$, which does not rule out variance fractions as large as $\Omega(n)$, though we keep the $1/n$ normalization so that $V_G$ is obtained from $E_G$ via the operator $\lambda \frac{\partial}{\partial \lambda}$.

Some early examples of extremal results on the hard-core model are bounds on the free energy. 
We state a theorem combining results of various well-known works.

\begin{theorem}[Cutler, Radcliffe, Kahn, Sah, Sawhney, Stoner, Zhao]\label{thm:Fbounds}
    Let $G$ be an $n$-vertex graph. Then for any $\lambda>0$
    \begin{equation}\label{eq:Ftrivial}
    \frac{1}{n}\log(1+n \lambda) = F_{K_n}(\lambda) \le F_{G}(\lambda) \le F_{\overline{K}_n}(\lambda) = \log(1+\lambda).
    \end{equation}
    If $G$ is $d$-regular then we also have
    \begin{equation}\label{eq:Fregular}
        \frac{1}{d+1}\log(1+(d+1)\lambda) = F_{K_{d+1}}(\lambda) \le F_{G}(\lambda) \le F_{K_{d,d}}(\lambda) = \frac{1}{2d}\log\left(2(1+\lambda)^d -1\right),
    \end{equation}
    and the lower bound holds under the weaker condition that $G$ has maximum degree $d$. 
    If $d_u$ is the degree of the vertex $u$ in $G$, then
    \begin{equation}\label{eq:Fdegrees}
        \frac{1}{n}\sum_{u\in V(G)} F_{K_{d_u+1}}(\lambda) \le F_{G}(\lambda) \le \frac{1}{n}\sum_{uv\in E(G)} \frac{d_u+d_v}{d_ud_v}F_{K_{d_u,d_v}}(\lambda) + \frac{1}{n}\sum_{\substack{u\in V(G) \\ \text{s.t. }d_u=0}}F_{K_1}(\lambda). 
    \end{equation}    
\end{theorem}

Note that~\eqref{eq:Ftrivial} is a trivial consequence of the simple fact that $F_G(\lambda) < F_{G-e}(\lambda)$ for any graph $G$ and any edge $e\in E(G)$, though we are not aware of it being stated explicitly in previous works. 
The lower bound in~\eqref{eq:Fregular} follows from a result of Cutler and Radcliffe~\cite{CR14} which we discuss in more detail later, and the upper bound combines celebrated results of Kahn~\cite{Kah01} and Zhao~\cite{Zha10}. 
The bounds in~\eqref{eq:Fdegrees} are a more recent development due to Sah, Sawhney, Stoner and Zhao~\cite{SSSZ19}. 

These results completely describe the extremal behavior of the free energy density in several important settings: fixing only the number of vertices, fixing the degree of every vertex to be $d$, and fixing an arbitrary degree sequence. 
In each case we have tight upper and lower bounds with explicit descriptions of the extremal graphs. 
The setting of bounded-degree graphs is also captured by the above bounds: the lower bound in~\eqref{eq:Fregular} holds under the condition of maximum degree $d$, showing that the complete graph $K_{d+1}$ minimizes the free energy density, and $K_1$ (equivalently, $\overline{K}_n$) is maximizer (the upper bound in~\eqref{eq:Ftrivial}).
From the degree sequence setting we also obtain sharp results for graphs of given average degree.

It is natural to consider generalizations of the above results to bounds on the parameters $E_G(\lambda)$ and $V_G(\lambda)$. 
Some such results are already known.

\begin{theorem}[Cutler, Radcliffe, Davies, Jenssen, Perkins, Roberts]\label{thm:Ebounds}
    Let $G$ be an $n$-vertex graph. Then for any $\lambda>0$
    \begin{equation}\label{eq:Etrivial}
    \frac{\lambda}{1+n\lambda} = E_{K_n}(\lambda) \le E_{G}(\lambda) \le E_{\overline{K}_n}(\lambda) = \frac{\lambda}{1+\lambda}.
    \end{equation}
    If $G$ is $d$-regular then we also have
    \begin{equation}\label{eq:Eregular}
        \frac{\lambda}{1+(d+1)\lambda} = E_{K_{d+1}}(\lambda) \le E_{G}(\lambda) \le E_{K_{d,d}}(\lambda) = \frac{\lambda(1+\lambda)^{d-1}}{2(1+\lambda)^d-1},
    \end{equation}
    and the lower bound holds under the weaker condition that $G$ has maximum degree $d$. 
\end{theorem}

The upper bound in~\eqref{eq:Etrivial} follows trivially from some conditional probability considerations, and the lower bound follows straightforwardly from results in the maximum degree setting.
The lower bound in~\eqref{eq:Eregular} can be proved in various ways, e.g.\ via the local occupancy methods of~\cite{DJPR17b} (see~\cite{Zha17,DK25}) or from the free volume arguments in~\cite{CR14}. We discuss these proofs later. 
Notably, bounds generalizing~\eqref{eq:Fdegrees} in the setting of a fixed degree sequence are absent from Theorem~\ref{thm:Ebounds}. 
This ``missing'' lower bound was recently conjectured by Davies and Kang~\cite[Conj.~A]{DK25}, but rather curiously the natural generalization of the upper bound in~\eqref{eq:Fdegrees} to the occupancy fraction is false, see Section~\ref{sec:FV}.

\subsection{Degree sequence bounds on occupancy fraction}

Our first result is a step towards the conjecture of Davies and Kang. While we obtain the desired bound, which is tight, we only do so for $\lambda$ sufficiently small. The full conjecture is that the same bound holds for all positive $\lambda$.

\begin{theorem}\label{thm:Edegrees}
    Let $G$ be a graph and let $d_u$ be the degree of a vertex $u$ in $G$. Then with $\Delta=\max_{u\in V(G)}d_u$, for any $\lambda \le \frac{3}{(\Delta+1)^2}$,
    \begin{equation}\label{eq:Edegrees}
        \frac{1}{n}\sum_{u\in V(G)} \frac{\lambda}{1+(d_u+1)\lambda} = \frac{1}{n}\sum_{u\in V(G)} E_{K_{d_u+1}}(\lambda) \le E_{G}(\lambda).
    \end{equation}
\end{theorem}

The Caro--Wei theorem~\cite{Caro79,Wei81} states that any graph $G$ contains an independent set of size at least 
\[ \sum_{u\in V(G)}\frac{1}{d_u+1}, \]
and a well-known proof due to Alon and Spencer~\cite{AS16} shows that a natural distribution achieves this in expectation: assign to each vertex $u$ independently, uniformly at random a value $x_u\in [0,1]$ and let $I = \{u\in V(G) : x_u < x_v \text{ for all } v\in N(u)\}$. 
In the limit $\lambda\to\infty$ the hard-core model approaches the uniform distribution over maximum independent sets, and thus~\eqref{eq:Edegrees} holds in the limit $\lambda\to\infty$ by the Caro--Wei theorem. 
That is, one can see Theorem~\ref{thm:Edegrees} as progress towards a smooth generalization of the Caro--Wei theorem for the hard-core model.
It is believable that a more careful analysis of our methods permits an upper bound of the form $\lambda = O(1/\Delta)$, but larger $\lambda$ may require new techniques. 

The proof of Theorem~\ref{thm:Edegrees} continues the development of a number of techniques and ideas due to Shearer~\cite{She83,She91}. The theme of these works is a lower bound on the independence number of a triangle-free graph in terms of its average degree and degree sequence respectively. 
The induction proof goes through if one can establish certain differential or difference equations subject to the local structure of the graph and the desired lower bound. 
Through the right lens, one can view certain proofs of the Caro--Wei theorem and related results in this way, though they are much easier. 
A very recent result of Buys, van den Heuvel and Kang~\cite{BvK25} generalizes Shearer's original method~\cite{She83} to the free energy density instead of the independence number (cf. PART vs MAX below). 
They push through a similar induction, but for the partition function instead of the independence number. 
Their proof is written for the restricted range $\lambda\in[0,1]$, but interestingly they are able to give an asymptotically tight lower bound on $F_G(1)$ in triangle-free graphs of average degree $d$.

It is natural to ask whether this approach can work for the occupancy fraction, and here we give an affirmative answer subject to a more significant restriction on the range of $\lambda$. 
This idea yields Theorem~\ref{thm:Edegrees} above, but we are also able to apply the method in the setting of triangle-free graphs and make some progress on~\cite[Conj.~12]{BvK25}. 
We show that a degree-sequence Shearer-like induction~\cite{She91} \emph{combined} with the local occupancy methods of~\cite{DJPR17a,DKPS20b} can give a lower bound on the occupancy fraction in terms of the degree sequence in the triangle-free setting. Since the function of $d_u$ appearing in the sum is convex, it gives an interesting result in the setting of given average degree too.

\begin{theorem}\label{thm:EdegreesTF}
    There is an absolute constant $c>0$ such that the following holds.
    Let $G$ be a triangle-free graph and let $d_u$ be the degree of a vertex $u$ in $G$. Then with $\Delta=\max_{u\in V(G)}d_u$, for any $\lambda \le c/\Delta^4$,
    \begin{equation}\label{eq:EdegreesTF}
        \frac{1}{n}\sum_{u\in V(G)} \frac{\lambda}{1+\lambda}\frac{ W(d_u \log(1+\lambda))}{d_u\log(1+\lambda)} \le E_{G}(\lambda),
    \end{equation}
    where $W$ is the (positive real branch of the) Lambert $W$-function.
\end{theorem}

This generalizes a bound of Davies, Jenssen, Perkins and Roberts~\cite{DJPR17a} which holds in the setting of fixed maximum degree, but here we have a limited range of $\lambda$.
For large $x$ it holds that $W(x)\sim\log x$, so the bound appears to be of the right order as $\Delta\to\infty$, but the condition on $\lambda$ precludes an application of this result in a range where the hard-core model reveals much about the structure of large independent sets in a graph.

Unlike Theorem~\ref{thm:Edegrees}, the bound in Theorem~\ref{thm:EdegreesTF} is not known to be tight in the stated range of $\lambda$. 
We find the question of how tight the bound is rather interesting, and given some analysis of the binomial random graph in~\cite{BvK25}, we expect that the lower bound we give is very close to what one finds in binomial random graph of constant average degree. In the case of regular graphs, tightness through comparison with the random regular graph was discussed in~\cite{DJPR17a}.
If we could significantly relax the upper bound on $\lambda$ to, say, $1/\log d$ where $d$ is the \emph{average} degree of $G$, then we would prove~\cite[Conj.~12]{BvK25}.

Our results on triangle-free graphs and for the variance (see below), as well as the methods developed in their pursuit, continue a line of work with applications to Ramsey theory. 
Shearer's induction method~\cite{She83} still provides the best-known upper bound on the Ramsey number $R(3,t)$, though one can also use local occupancy~\cite{DJPR17a} to get a bound of the same leading order asymptotically. 
The best lower bound is a factor $3+o(1)$ smaller and comes from an analysis of an elegant construction that balances structure and randomness carefully~\cite{CJMS25}, improving upon famous prior work analyzing the triangle-free process~\cite{BK21a,FGM20}.
A conjecture in~\cite{DJPR17a} suggests that through variance bounds on the hard-core model one could hope to improve the leading order of the upper bound. 
The gap in the bounds on $R(3,t)$ is one of the most well-known open problems in Ramsey theory, and the fact that Shearer's bound remains unsurpassed despite decades of sustained effort is evidence that the program we pursue here is fraught with difficulty. 
We find any prospect of connections between such gaps and the aforementioned algorithmic barrier interesting too.

\subsection{Bounds on variance fraction}

Our third result concerns generalizations of Theorem~\ref{thm:Fbounds} to the variance fraction. 
Here, much less is known. One of the essential difficulties is that the variance is the centered second moment, and by Theorems~\ref{thm:Ebounds} and~\ref{thm:Edegrees} the very graphs which we hope to show minimize the variance also minimize the magnitude of the negative term corresponding to the centering. 
By formulating a higher-order version of the occupancy method, we establish a first step towards the simplest general bounds on the variance fraction. 
As in the preceding subsection, our results are limited to small values of $\lambda$.

\begin{theorem}\label{thm:Vbounds}
    Let $G$ be an $n$-vertex graph. Then for any $0<\lambda<1/(2n-1)$ we have
    \begin{equation}\label{eq:Vtrivial}
    \frac{\lambda}{(1+n\lambda)^2} = V_{K_n}(\lambda) \le V_{G}(\lambda) \le V_{\overline{K}_n}(\lambda) = \frac{\lambda}{(1+\lambda)^2},
    \end{equation}
    and the upper bound holds up to $\lambda \le 1/n$.
\end{theorem}

We conjecture an extension of the lower bound to arbitrary positive $\lambda$ and warn that the upper bound cannot hold for all graphs and all $\lambda>0$, see Section~\ref{sec:VAR}.

\begin{conjecture}\label{conj:V}
    Let $G$ be an $n$-vertex graph. Then for any $\lambda>0$ we have
    \begin{equation}\label{eq:Vtrivialconj}
    \frac{\lambda}{(1+n\lambda)^2} = V_{K_n}(\lambda) \le V_{G}(\lambda)
    \end{equation}
    If $G$ has maximum degree $\Delta$ then the lower bound can be improved to 
    \begin{equation}\label{eq:Vregularconj}
    \frac{\lambda}{(1+(\Delta+1)\lambda)^2} = V_{K_{\Delta+1}}(\lambda) \le V_{G}(\lambda).
    \end{equation}
\end{conjecture}

It is plausible that the generalization of~\eqref{eq:Edegrees} holds for variance too. 
We discuss the smallest nontrivial case in Section~\ref{sec:FV}.

\subsection{Free volume and related inequalities}

The partition function $Z_G(x)$ is a monic polynomial in $x$ with nonnegative (integer) coefficients. 
There are a variety of extremal results that one might seek in this context, and a taxonomy of such results and their relationships was studied in~\cite{DJPR18} (among others). 

For a fixed $n$, consider polynomials $P(x) = \sum_{k=0}^n a_k x^k$ and $Q(x) = \sum_{k=0}^n b_k x^k$. 
Suppose that $a_0 = b_0 = 1$ and $a_k, b_k \ge 0$ for all $k$. 
Note that we allow coefficients to be zero, so $P$ and $Q$ may not technically have the same degree. 
This is important when considering our primary application to the hard-core model: two comparable graphs, e.g.\ on the same number of vertices, may not have the same independence number.
The following statements were studied in~\cite{DJPR18}.

\begin{definition}\hfill
    \begin{enumerate}
        \item We say that $P\geCOUNT Q$ if $ \sum_{k=0}^n a_k \ge \sum_{k=0}^n b_k$.
        \item We say that $P\gePART Q$ if $ P(x) \ge Q(x)$ for all $x \ge 0$.
        \item We say that $P\geCOEF Q$ if $ a_k \ge b_k$ for all $1 \le k \le n$.
        \item We say that $P\geOCC Q$ if $ \frac{x P'(x)}{P(x)} \ge \frac{x Q'(x)}{Q(x)}$ for all~$x \ge 0$.
        \item We say that $P\geMAX Q$ if $ a_n \ge b_n$.
        \item We say that $P\geFV Q$ if $ b_ka_{k+1}\ge a_kb_{k+1}$ for all $0 \le k \le n-1$.     
    \end{enumerate}
\end{definition}

These definitions capture a large number of results and conjectures from the literature on independent sets, matchings, and other combinatorial structures. 
For example, Bregman's theorem~\cite{Bre73} is a MAX-type bound on the partition function of the monomer-dimer model, and the upper matching conjecture~\cite{FKM08} is the conjectured extension to COEF in the $d$-regular setting. 
See~\cite{DJP21} for a proof of this under the additional assumption that the graph has $n$ vertices for $n\ge n_0(d)$ that makes use of OCC-type bounds from~\cite{DJPR17b}.

Almost all of the results mentioned so far fit into this framework, especially if we add an extension to cover variance.
Given our definitions of free energy and occupancy fraction we have that $Z_G \gePART Z_H$ is equivalent to $F_G(x) \ge F_H(x)$ for all $x \ge 0$ and $Z_G \geOCC Z_H$ is equivalent to $E_G(x)\ge E_H(x)$ for all $x \ge 0$. 
Thus, it is natural to add the statement corresponding to $V_G(x)\ge V_H(x)$ for all $x \ge 0$.
\begin{definition}
Given a twice differentiable function $P:\mathbb{R}\to\mathbb{R}$ we write \[V_P(x) = \frac{x^2 P''(x) + x P'(x)}{P(x)} - \frac{x^2 P'(x)^2}{P(x)^2}.\] 
To the list of properties above, we add the following.
    \begin{enumerate}\foreach \n in {1,...,6} {\refstepcounter{enumi}}
        \item We say that $P\geVAR Q$ if $V_P(x) \ge V_Q(x)$ for all $x\ge 0$.
    \end{enumerate}
\end{definition}
We relate VAR to the other properties as follows.

\begin{theorem}
    For a fixed $n$, let $P(x) = \sum_{k=0}^n a_k x^k$ and $Q(x) = \sum_{k=0}^n b_k x^k$ be polynomials. 
    Suppose that $a_0 = b_0 = 1$ and $a_k, b_k \ge 0$ for all $k$. Then 
    \begin{enumerate}
        \item $P \geVAR Q \Rightarrow P \geOCC Q$.
        \item FV and VAR are incomparable in general in the sense that there exists a choice of polynomials as above such that $P \geFV Q$ but $ P \not\geVAR Q$, and there exists a different choice such that $P \not\geFV Q$ but $ P \geVAR Q$.
        \item COEF and VAR are incomparable in general in the sense that there exists a choice of polynomials as above such that $P \geCOEF Q$ but $ P \not\geVAR Q$, and there exists a different choice such that $P \not\geCOEF Q$ but $ P \geVAR Q$.
    \end{enumerate}
\end{theorem}

\noindent
In conjunction with~\cite[Prop.~19]{DJPR18}, we have the following web of implications.

\begin{figure}[H]
    \centering
    \begin{tikzpicture}[scale=1.25]
        \tikzset{sn/.style={rectangle,inner sep=2pt,minimum width=50pt,minimum height=15pt}}
        \node[sn] (var) at (2,0) {VAR};
        \node[sn] (fv) at (0,0) {FV};
        \node[sn] (coef) at (-1,-1) {COEF};
        \node[sn] (occ) at (1,-1) {OCC};
        \node[sn] (part) at (0,-2) {PART};
        \node[sn] (count) at (-1,-3) {COUNT};
        \node[sn] (max) at (1,-3) {MAX};

        \path (fv) -- (var) node[midway,sloped] {$\substack{\not\Rightarrow\\\not\Leftarrow}$};
        \path (var) -- (occ) node[midway,sloped] {$\Leftarrow$};
        \path (fv) -- (coef) node[midway,sloped] {$\Leftarrow$};
        \path (fv) -- (occ) node[midway,sloped] {$\Rightarrow$};
        \path (coef) -- (occ) node[midway,sloped] {$\substack{\not\Rightarrow\\\not\Leftarrow}$};
        \path (coef) -- (part) node[midway,sloped] {$\Rightarrow$};
        \path (occ) -- (part) node[midway,sloped] {$\Leftarrow$};
        \path (part) -- (count) node[midway,sloped] {$\Leftarrow$};
        \path (part) -- (max) node[midway,sloped] {$\Rightarrow$};
        \path (count) -- (max) node[midway,sloped] {$\substack{\not\Rightarrow\\\not\Leftarrow}$};
    \end{tikzpicture}
    \caption{Implications between extremal results for polynomials with nonnegative coefficients.\label{fig:implications}}
\end{figure}
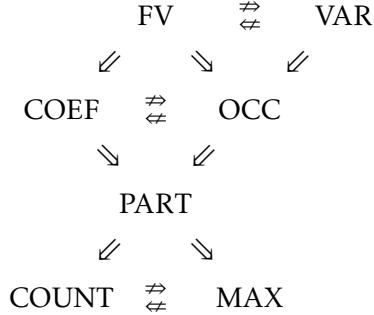

It is plausible that a generalization of Theorem~\ref{thm:Edegrees} holds in the sense that the free volume inequalities $P \geFV Q$ hold when $P=Z_G(\lambda)$ and $Q=\prod_{u\in V(G)}Z_{K_{d_u+1}}(\lambda)^{\frac{1}{d_u+1}}$, though one has to take care with this statement as $Q$ is only a polynomial if for each $d$ the number of vertices of degree $d$ in $G$ is divisible by $d+1$.
We verified the first nontrivial case of this type of bound in Section~\ref{sec:FV}.

\subsection{Bounds combining parameters}

We conclude the introduction with a reference to an interesting result of Campos and Samotij~\cite{CS24a}. 
This serves to highlight preexisting interest in extremal results which relate the parameters discussed thus far. 
We are aware of two results relating $E_G$ to $F_G$. 

\begin{theorem}[Campos, Samotij, Shearer, Davies, Kang, Pirot, Sereni]
    Let $G$ be a graph on $n$ vertices and $x>0$. Then 
    \begin{equation}\label{eq:combined}
        \frac{(1+\lambda)\log(1+\lambda)}{\lambda} E_G(\lambda) \le F_G(\lambda) \le E_G(\lambda)\log(\lambda) + h\left(E_G(\lambda)\right) \le E_G(\lambda)\log\left(e\lambda/E_G(\lambda)\right),
    \end{equation}
    where $h:[0,1]\to\mathbb{R}$ is the base-$e$ entropy function given by $h(0)=h(1)=0$ and $h(x) = -x\log x - (1-x)\log(1-x)$.
\end{theorem}
The first inequality is due to Campos and Samotij~\cite{CS24a} and is tight for a graph with no edges. The second is a generalization of an argument of Shearer~\cite{She95} first published in~\cite{DKPS20b}, and we do not expect it to be tight on any finite graph.

The first inequality in~\eqref{eq:combined} strengthens the upper bound on $E_G(\lambda)$ in~\eqref{eq:Etrivial} because by~\eqref{eq:Ftrivial} we have $F_G(\lambda)\le \log(1+\lambda)$.
Given this, one might be tempted to investigate the derivative of the first inequality. This is equivalent to $(1+\lambda)V_G(\lambda) \le E_G(\lambda)$, but such a bound is false on the smallest graph which is not a disjoint union of cliques: $K_{1,2}$ (at any $\lambda>1+\sqrt 5$).

\section{Degree sequence bounds on occupancy fraction}\label{sec:degrees}

We begin this section with a discussion of local occupancy~\cite{DJPR17a,DKPS20a,DK25} as it is an essential tool in the proofs of Theorems~\ref{thm:Edegrees} and~\ref{thm:EdegreesTF}.

Given a graph $G$, and real parameters $\lambda,\beta,\gamma>0$, we say that $G$ has \emph{local $(\beta,\gamma)$-occupancy at fugacity $\lambda$} if for each $u\in V(G)$ and induced subgraph $F\subset G[N(u)]$ we have 
\begin{equation}\label{eq:lococc} \beta\frac{\lambda}{1+\lambda}\frac{1}{Z_F(\lambda)} + \gamma \frac{\lambda Z_F'(\lambda)}{Z_F(\lambda)}. \end{equation}

It has been established before~\cite{DJPR17a,DKPS20a} that in graphs of maximum degree $\Delta$ which satisfy local  $(\beta,\gamma)$-occupancy at fugacity $\lambda$ we have 
\[ E_G(\lambda) \ge \frac{1}{\beta+\Delta\gamma},\]
and that every graph satisfies local $(1+1/\lambda, 1)$-occupancy at fugacity $\lambda$ for every $\lambda>0$~\cite{DK25}. 
Note that this gives a proof of the lower bounds in~\eqref{eq:Etrivial} and~\eqref{eq:Eregular} (and thus also~\eqref{eq:Ftrivial} and~\eqref{eq:Fregular} via some integration).
What we require here is the following slightly stronger bound that follows from the same approach. 

\begin{theorem}\label{thm:lococclb}
    Suppose that $G$ satisfies local $(\beta,\gamma)$-occupancy at fugacity $\lambda$. Then writing $I$ for a random independent set from the hard-core model $\mu_{G,\lambda}$ on $G$ at fugacity $\lambda$,
    \begin{align*}
        \frac{1}{n}\sum_{u\in V(G)}\Pr(u\in I)\cdot \left(\beta+d_u\gamma\right) \ge 1.
    \end{align*}
\end{theorem}
\begin{proof}
    The spatial Markov property of the model combined with~\eqref{eq:lococc} give the result.
    For each $u$, we consider the random subgraph $F_u$ of $G[N(u)]$ obtained as follows. 
    Sample $I\sim \mu_{G,\lambda}$ and reveal $I\setminus N(u)$. Then let $F_u$ be the subgraph of $G[N(u)]$ induced by the set 
    $\{ v\in N(u) : N(v)\cap (I\setminus N(u)) = \emptyset \}$
    of vertices uncovered by $I\setminus N(u)$.

    Standard computations in the hard-core model (e.g.~\cite{DJPR17b,DK25}) give that 
    \begin{align*} \Pr(u\in I \mid F_u=F) &= \frac{\lambda}{1+\lambda}\frac{1}{Z_{F}(\lambda)},&
    \sum_{v\in N(u)}\Pr(v\in I \mid F_u=F) &= \frac{\lambda Z_{F}'(\lambda)}{Z_{F}(\lambda)}. \end{align*}
    But combining these bounds with the law of total expectation and~\eqref{eq:lococc}, we have for every $u$,
    \[ \beta \Pr(u\in I) + \gamma \sum_{v\in N(u)} \Pr(v\in I) \ge 1. \]
    Averaging this inequality over a uniform random choice of vertex $u$ gives the result, after interchanging the order of the sums to see that
    \[ \sum_{u\in V(G)}\sum_{v\in N(u)}\Pr(v\in I) = \sum_{u\in  V(G)}d_u \Pr(u\in I).\qedhere\]
\end{proof}

Combined with the aforementioned local occupancy parameters from~\cite{DK25}, we have the following result.

\begin{corollary}\label{cor:genlococc}
    Let $G$ be an $n$-vertex graph and $\lambda>0$. Let $f_\lambda(d)$ be the function $Z_{K_{d+1}}(\lambda)=\frac{\lambda}{1+(d+1)\lambda}$. Then if $I$ is a random independent set from the hard-core model on $G$ at fugacity $\lambda$,
    \[ \frac{1}{n}\sum_{u\in V(G)}\frac{\Pr(u\in I)}{f_\lambda(d_u)} \ge 1. \]
\end{corollary}
\begin{proof}
    The result follows from Theorem~\ref{thm:lococclb} and the fact~\cite{DK25} that every graph satisfies local $(1+1/\lambda, 1)$-occupancy at fugacity $\lambda$ for every $\lambda>0$. The latter is an easy consequence of Markov's inequality.
\end{proof}

This resembles in a superficial way the desired bound in the degree-sequence setting. Theorem~\ref{thm:Edegrees} is the bound
\[ E_G(\lambda) = \frac{1}{n}\sum_{u\in V(G)}\Pr(u\in I) \ge \frac{1}{n}\sum_{u\in V(G)} f_\lambda(d_u), \]
but we have some work to do to untangle the marginals from the function $f_\lambda$ in the average.
It is in performing the untangling that we require $\lambda$ to be small.

\subsection{General graphs: proof of Theorem~\ref{thm:Edegrees}}\label{sec:degreesGEN}

\begin{proof}[Proof of Theorem~\ref{thm:Edegrees}]
    It suffices to consider connected graphs because the occupancy fraction of a disjoint union of graphs is a sum of the individual occupancy fractions weighted appropriately by their order. In particular, we may assume that $G$ is either $K_1$ or has minimum degree $1$.
    
    Let $\Delta$ be the maximum degree of $G$ and $V=V(G)$. We proceed by induction on $n=|V|$. Note that subgraphs of $G$ cannot have higher maximum degree that $G$, and hence if we start with $G$ and $\lambda$ such that $\lambda \le 3/(\Delta+1)^2$, then the analogous condition holds for all subgraphs of $G$ that we see in the induction. 
    For the entire proof $\lambda$ remains fixed and so we drop it from the notation, writing $f(d) = f_{\lambda}(d)$, $\mu_{G}=\mu_{G,\lambda}$, etc. The base case is a graph of one vertex in which there is nothing to prove: the occupancy fraction is $f(0)$ as required. 

    For the induction, we consider an arbitrary vertex $u\in V(G)$ and compute the occupancy fraction conditioned on the event that $u$ is occupied. 
    Let $p_v$ be the probability that a vertex $v$ is occupied when an independent set in $G$ is drawn from the hard-core measure $\mu_{G}$. 
    By the spatial Markov property of the measure, conditioning on $u\notin I$ results in a set distributed according to $\mu_{G-u}$ and conditioning on $u\in I$ we get $I\setminus \{u\}$ distributed according to $\mu_{G-N[u]}$. 
    Then with $\E_{G}$ meaning an expectation with respect to $\mu_{G}$ and $d_{uw}=|N(u)\cap N(w)|$ (in the graph $G$) we obtain the following by induction. 
    The variables $v,w,x,yz$ in the sums below always range over $v\in N(u)$, $w\in N^2(u)$, $x\in V\setminus N[u]$, $y\in V\setminus N^2[u]$, and $z\in V\setminus\{u\}$. 
    \begin{align*}
    \E_{G}|I|  &= (1-p_u)\E_{G-u}|I| + p_u\left(1+\E_{G-N[u]}|I|\right)
    \\&\ge (1-p_u)\left(\sum_{v}f(d_v-1) + \sum_{x}f(d_x)\right) + p_u\left(1 + \sum_{w}f(d_w-d_{uw}) + \sum_{y}f(d_y)\right)
    \\&= \sum_{z}f(d_z) + p_u\left(1+\sum_{w}[f(d_w-d_{uw})-f(d_w)] - \sum_{v}f(d_v-1)\right) + \sum_{v}[f(d_v-1)-f(d_v)].
    \end{align*} 
Thus, we are done if there exists a vertex $u$ such that 
\[ p_u\left(1+\sum_{w}[f(d_w-d_{uw})-f(d_w)] - \sum_{v}f(d_v-1)\right) \ge  f(d_u) - \sum_{v}[f(d_v-1)-f(d_v)]. \]

Since $f$ is decreasing it suffices to show that there exists a vertex $u$ such that 
\[ p_u\left(1- \sum_{v\in N(u)}f(d_v-1)\right) \ge  f(d_u) - \sum_{v\in N(u)}[f(d_v-1)-f(d_v)]. \]
Since $d_v\ge 1$ in the sums above, $f$ is decreasing, and $f(0)=\lambda/(1+\lambda)$, our upper bound on $\lambda$ gives that the coefficient of $p_u$ is strictly positive. Then we can divide by it and observe that it suffices to show that there exists a vertex $u$ such that 
\[ \frac{p_u}{f(d_u)} \ge  \frac{1 - \sum_{v\in N(u)}\frac{f(d_v-1)-f(d_v)}{f(d_u)}}{1- \sum_{v\in N(u)}f(d_v-1)}. \]

As is standard in Shearer-type inductions, we show that the above inequality holds on average over a uniform random choice of $u$. 
By Corollary~\ref{cor:genlococc}, the average of the left-hand side is at least 1, so it suffices to show that 
\begin{equation}\label{eq:precrazyfact}
    1 \ge \frac{1}{n}\sum_{u\in V(G)}\frac{1 - \sum_{v\in N(u)}\frac{f(d_v-1)-f(d_v)}{f(d_u)}}{1- \sum_{v\in N(u)}f(d_v-1)}.
\end{equation}
The use of local occupancy to average out the probabilities that appear when one conditions on $u\in I$ and $u\notin I$ is novel, this difficulty does not arise in analogous proofs involving independence number or free energy density.

We proceed via the surprising identity 
\begin{equation}\label{eq:fidentity}\frac{f(d_v-1)-f(d_v)}{f(d_u)} = f(d_v-1) + (d_u-d_v)f(d_v-1)f(d_v), \end{equation}
which holds for all $d_u$, $d_v$, and $\lambda$. 
Using this in~\eqref{eq:precrazyfact} and canceling terms, it suffices to show that
\begin{equation*}\
    0 \le \frac{1}{n}\sum_{u\in V(G)}\sum_{v\in N(u)}\frac{(d_u-d_v)f(d_v-1)f(d_v)}{1- \sum_{v\in N(u)}f(d_v-1)}.
\end{equation*}
But we can write the double sum as a sum over edges instead to get the equivalent expression
\begin{align*}
    \sum_{\substack{uv\in E \\ d_u > d_v}}(d_u-d_v)\left( \frac{f(d_v-1)f(d_v)}{1-\sum_{v'\in N(u)}f(d_{v'}-1)} - \frac{f(d_u-1)f(d_u)}{1-\sum_{u'\in N(v)}f(d_{u'}-1)}\right) \ge 0.
\end{align*}
Note that edges between vertices of the same degree contribute nothing, so we can sum over the edges $uv$ such that $d_u-d_v> 0$.
Since $f$ is strictly positive we get a lower bound on the left-hand side by replacing the denominator
\[1-\sum_{v'\in N(u)}f(d_{v'}-1) \]
with the upper bound
\[ 1 - d_u\lambda + d_v\lambda^2 + \sum_{v'\in N(u)\setminus\{v\}}d_{v'}\lambda^2 \le 1 - d_u\lambda + d_v\lambda^2 + (d_u-1)\Delta\lambda^2 \]
that follows from the simple fact that $f(d-1) \ge \lambda - d\lambda^2$ for all $d,\lambda\ge0$.
Similarly, we get a lower bound on the left-hand side by replacing the other denominator $1-\sum_{u'\in N(v)}f(d_{u'}-1)$.
with the lower bound $1- d_v \lambda$.
Then we are done if 
\begin{align*}
    \sum_{\substack{uv\in E \\ d_u > d_v}}(d_u-d_v)\left( \frac{f(d_v-1)f(d_v)}{1 - d_u\lambda + d_v\lambda^2 + (d_u-1)\Delta\lambda^2} - \frac{f(d_u-1)f(d_u)}{1- d_v \lambda}\right) \ge 0.
\end{align*}
Since $f$ is decreasing, each term in the summand is a decreasing function of $d_v$. Then we get a lower bound on each term by replacing $d_v$ with the maximum value $d_u-1$.
Then it suffices to show that
\begin{align*}
    \sum_{\substack{uv\in E \\ d_u > d_v}}(d_u-d_v)f(d_u-1)\left( \frac{f(d_u-2)}{1 - d_u\lambda + (d_u-1)(\Delta+1)\lambda^2} - \frac{f(d_u)}{1- (d_u-1) \lambda}\right) \ge 0.
\end{align*}
Note that $d_u\ge2$ for all edges $uv$ in the sum.
Recall that our lower bound on $\lambda$ means that the denominators are positive, so we merely need for every $d_u\ge 2$ that 
\begin{align*}
f(d_u-2)\left(1-(d_u-1)\lambda\right) &\ge f(d_u)\left(1-d_u\lambda +(d_u-1)(\Delta+1)\lambda^2\right) &&\Leftrightarrow\\
\left(1+(d_u+1)\lambda\right)\left(1-(d_u-1)\lambda\right) &\ge \left(1+(d_u-1)\lambda\right)\left(1-d_u\lambda +(d_u-1)(\Delta+1)\lambda^2\right)&& \Leftrightarrow \\
1 + 2\lambda - (d_u^2-1)\lambda^2 &\ge 1 - \lambda + (\Delta+1 - d_u)(d_u-1)\lambda^2&& \Leftrightarrow \\
3\lambda &\ge (\Delta+2)(d_u-1)\lambda^2.
\end{align*}
This is implied by the upper bound $\lambda \le 3/(\Delta+1)^2$, as required.
\end{proof}

\subsection{Triangle-free graphs: proof of Theorem~\ref{thm:EdegreesTF}}\label{sec:degreesTF}

The proof for the triangle-free case is essentially the same, but one has to work a lot harder as the functions involved are more difficult to expand in $\lambda$. One also has to exploit triangle-freeness, which comes in two forms. 
First, we get to apply a stronger lower bound derived from a local occupancy analysis that exploits triangle-freeness, and second when we tame the resulting sums we apply the following fact. In any triangle-free graph, we have for all vertices $u$ that
\begin{equation}\label{eq:tfedgecounts} 
\sum_{v\in N(u)} d_v  = d_u + \sum_{w\in N^2(u)} d_{uw}.
\end{equation}
This follows from the fact that $N(u)$ is an independent set in $G$ and so all edges incident to a neighbor of $u$ must also be incident to a vertex in the set $\{u\}\cup N(w)$.

For this subsection, let 
\begin{align*} g(d) &= \frac{\lambda}{1+\lambda} \frac{W(d \log(1+\lambda))}{d\log(1+\lambda)},
\end{align*} 
where $W$ is the Lambert $W$-function.
This $W$ is the inverse of $x\mapsto x e^x$, and we need the branch defined on $[0,\infty)$,
It was shown in~\cite{DJPR17a} that for every $\lambda> 0$, every triangle-free graph satisfies local $(\beta,\gamma)$-occupancy with parameters such that $\beta + d_u \gamma = g(d_u)$.  
We state the necessary corollary of this fact and Theorem~\ref{thm:lococclb} here.

\begin{corollary}\label{cor:TFlococc}
    Let $G$ be an $n$-vertex triangle-free graph and $\lambda>0$. Let $g(d)$ be the function above. Then if $I$ is a random independent set from the hard-core model on $G$ at fugacity $\lambda$,
    \[ \frac{1}{n}\sum_{u\in V(G)}\frac{\Pr(u\in I)}{g(d_u)} \ge 1. \]
\end{corollary}

Since we rely on local occupancy in the proof this choice of $g$ is essentially forced upon us, though we could choose any function upper bounded by $g$. 
It would be interesting to compare our function with that in~\cite{BvK25}, where there appears to be more flexibility in the method. 
Given a function such as $g$, an identity analogous to~\eqref{eq:fidentity} seems unlikely to come to our aid. 
Circumventing this seems to contribute to the worsening of the upper bound on $\lambda$ in Theorem~\ref{thm:EdegreesTF} compared to Theorem~\ref{thm:Edegrees}.
Analogous to the improvement Shearer made when sharpening his result in terms of average degree~\cite{She83} to one in terms of degree sequence~\cite{She91}, we speculate that there is a lower bound on the occupancy fraction of triangle-free graphs that is more suited to difference equations than $g$.

\begin{proof}[Proof of Theorem~\ref{thm:EdegreesTF}]
    It suffices to consider connected graphs because the occupancy fraction of a disjoint union of graphs is a sum of the individual occupancy fractions weighted appropriately by their order. In particular, we may assume that $G$ is either $K_1$ or has minimum degree $1$.
    
The method is a similar induction to the previous proof. The base case is $n=1$, in which case the result holds by the definition of $g(0)$.

For the induction, we are done if there exists a vertex $u$ such that
\[ p_u\left(1+\sum_{w\in N^2(u)}[g(d_w-d_{uw})-g(d_w)] - \sum_{v\in N(u)}g(d_v-1)\right) \ge  g(d_u) - \sum_{v\in N(u)}[g(d_v-1)-g(d_v)]. \]
We restrict $\lambda$ such that the coefficient of $p_u$ is strictly positive. Since $f$ is decreasing, the coefficient is at least $1 - d_u g(0) \ge 1 - \Delta \lambda$, and hence $\lambda < 1/\Delta$ suffices for this step.

Then, after some algebraic manipulation, averaging over a uniform random vertex $u$, and after an application of Corollary~\ref{cor:TFlococc}, it suffices to show that 
\begin{align} \label{eq:tfwant}
1 &\ge \frac{1}{n}\sum_{u\in V} \frac{1 - \sum_{v\in N(u)}\frac{g(d_v-1)-g(d_v)}{g(d_u)}}{1+\sum_{w\in N^2(u)}[g(d_w-d_{uw})-g(d_w)] - \sum_{v\in N(u)}g(d_v-1)}.
\end{align}
We establish a stronger inequality, obtained as follows. 

First, we replace the term $t=\left(g(d_v-1)-g(d_v)\right)/{g(d_u)}$ in the numerator with a lower bound. For this, we compute a Taylor expansion in $\lambda$ at zero.
It is easy to prove that the coefficient of $\lambda^k$ in the expansion of $t$ is a polynomial in $d_u$ and $d_v$ of degree at most $k-1$. Then we have 
\[ t = a_0 + a_1\lambda + a_2\lambda^2 + a_3\lambda^3 + a_4 \lambda^4 + O(\Delta^4\lambda^5). \]
Then, there exists an absolute constant $c$ such that we get a lower bound on $t$ valid for all $0<\lambda < c/\Delta$ by truncating the series at the term $\lambda^4$ and replacing the true coefficient $a_4$ of $\lambda^4$ by one that is at most $a_4-\Delta^3$. 
We compute\begin{align*}
12 a_4 &= 18 d_u^2 d_v+96 d_u d_v^2-42 d_u d_v+8 d_u^3+16 d_u-250 d_v^3+231 d_v^2-139d_v+28 \\
&\ge -42 d_u d_v-250 d_v^3-139 d_v \ge - 431 \Delta^3,
\end{align*}
and using that $431/12+1 < 37$ we have that 
\[ \tfrac{g(d_v-1)-g(d_v)}{g(d_u)} \ge \lambda + \left(d_u+1-3d_v\right)\lambda^2 + \frac{1}{2}\left(3+d_u-d_u^2-10d_v-6d_ud_v+16d_v^2\right)\lambda^3 - 37\Delta^3\lambda^4. \]
Since the left-hand side is positive, we also restrict $\lambda$ such that our lower bound on the right is positive. 
Crudely, it is easy to show that this holds for all $\lambda \le 1/(50\Delta)$.

Second, we replace the term $t'=g(d_w-d_{uw})-g(d_w)$ in the denominator with a lower bound. Again, we proceed with a Taylor expansion in $\lambda$ at zero.
It is easy to prove that the coefficient of $\lambda^k$ in the expansion of $t'$ is a polynomial in $d_w$ and $d_{uw}$ of degree at most $k-1$. 
Writing
\[ t' = a'_0 + a'_1\lambda + a'_2\lambda^2 + a'_3\lambda^3 + a'_4 \lambda^4 + O(\Delta^4\lambda^5), \]
we have 
\[ a'_4 = \frac{8 d_{u w}^3}{3}-8 d_w d_{u w}^2-3 d_{u w}^2+8 d_w^2 d_{u w}+6 d_w d_{u w}+\frac{11 d_{u w}}{6}. \]
It is straightforward to see that $a'_4$ is an increasing function of $d_{uw}$ (which must be between $1$ and $d_{w}$) and is therefore minimized at the lower limit $d_{uw}=1$:
\[ \frac{\partial}{\partial d_{uw}}a'_4 = \frac{11}{6} + 8(d_w-d_{uw})^2 + 6(d_w-d_{uw}). \]
We then compute that $a'_4 \ge 11/8$ for all $d_w\ge 0$.
Then, for some absolute constant $c$ we get a lower bound on $t'$ valid for $\lambda \le c/\Delta$ by truncating the expansion at the $\lambda^3$ term. That is,
\[ t' \ge \lambda ^2 d_{u w}-\frac{3}{2} \lambda ^3 \left(-d_{u w}^2+2 d_w d_{u w}+d_{uw}\right) \ge d_{uw}(\lambda^2-3d_w\lambda^3), \]
where we use that $d_{uw}\ge1$ to get rid of the term $d_{uw}^2$ that would be hard to handle later in the proof. 
This lower bound is positive when $\lambda\le 1/(4\Delta)$.

Third, we replace the term $t''=g(d_v-1)$ in the denominator with an upper bound. 
By an argument very similar to the above, we get an upper bound valid in a range $\lambda \le O(1/\Delta)$ by truncating the Taylor expansion at the $\lambda^3$ term.
This gives  
\[ t'' \le \lambda - d_v \lambda ^2 + \frac{1}{2} \left(3 d_v^2-3 d_v+2\right) \lambda ^3. \]
Note that this is a stronger upper bound than $\lambda$ for all $\Delta\ge1$ and $\lambda\le 1/(4\Delta)$. 

After making these replacements, we have an upper bound on the right-hand side of~\eqref{eq:tfwant} that shares the same Taylor expansion with the true right-hand side up to and including the term quadratic in $\lambda$. 
We lost some precision in the $\lambda^3$ term by replacing some $d_{uw}^2$ terms with $d_{uw}$, but we still have sufficient accuracy here. Note that this replacement is tight whenever $d_{uw}=1$ which, for example, holds for all vertices $u$ and $w\in N^2(u)$ in any graph of girth at least five.

The rest of the proof is conceptually straightforward, if arduous. We keep performing Taylor-like expansions in such a way that we have perfect fidelity in terms up to $\lambda^3$.

We want to show that
\begin{align*}
    1 &\ge \frac{1}{n}\sum_{u\in V} \tfrac{1 - \sum_{v\in N(u)}\left(\lambda + \left(d_u+1-3d_v\right)\lambda^2 + \frac{1}{2}\left(3+d_u-d_u^2-10d_v-6d_ud_v+16d_v^2\right)\lambda^3 - 36\Delta^3\lambda^4\right)}{1 + \sum_{w\in N^2(u)} d_{uw}\left(\lambda^2-3d_w\lambda^3\right) - \sum_{v\in N(u)}\left(\lambda - d_v\lambda^2 + \frac{1}{2} (3d_v^2-3d_v+2)\lambda^3\right)}.
\end{align*}
To handle the denominator in a way that is accurate to terms of order $\lambda^3$,  we use that for all $y$ with $|y|\le 1/2$ we have 
\[ 1 + y + y^2 + y^3 + 2y^4 + 2y^5 \ge \frac{1}{1-y}.\]
Then it suffices to show that
\begin{align}\label{eq:TFwantatlast}
    1 &\ge \frac{1}{n}\sum_u\left(1 - x_u\right)\left(1 + y_u + y_u^2 + y_u^3 + 2y_u^4+ 2y_u^5\right),
\end{align}
where 
\begin{align*}
x_u &=\sum_v\left(\lambda + \left(d_u+1-3d_v\right)\lambda^2 + \frac{1}{2}\left(3+d_u-d_u^2-10d_v-6d_ud_v+16d_v^2\right)\lambda^3 - 37\Delta^3\lambda^4\right)
\\&= d_u\lambda + \sum_v\left(\left(d_u+1-3d_v\right)\lambda^2 + \frac{1}{2}\left(3+d_u-d_u^2-10d_v-6d_ud_v+16d_v^2\right)\lambda^3 - 37\Delta^3\lambda^4\right),
\shortintertext{and}
y_u &= -\sum_w d_{uw}\left(\lambda^2-3d_w\lambda^3\right) + \sum_v\left(\lambda - d_v\lambda^2 + \frac{1}{2} (3d_v^2-3d_v+2)\lambda^3\right)
\\&=d_u\lambda + \left(d_u-2\sum_v d_v\right)\lambda^2 + \left(\sum_v\frac{1}{2} (3d_v^2-3d_v+2) + 3\sum_w d_wd_{uw}\right)\lambda^3,
\end{align*}
provided we also show that $|y_u|\le 1/2$.
Note that in the equality for $y_u$ above, we use~\eqref{eq:tfedgecounts}.

The right-hand side of~\eqref{eq:TFwantatlast} is a polynomial in $\lambda$ and can be written as
\[ b_0 + b_1\lambda + b_2\lambda^2 + b_3\lambda^3 + \xi, \]
where the $b_i$ are independent of $\lambda$ and $\lambda^4$ divides the polynomial $\xi$ in $\lambda$.
In what follows, we abbreviate the summation notation for convenience. In sums, we always have $u\in V$, $v\in N(u)$, and $w\in N^2(u)$.
We compute
\begin{align*}
    b_0 &= 1\\
    b_1 &= \frac{1}{n}\sum_u\left(- d_u + d_u\right) = 0,\\
    b_2 &= \frac{1}{n}\sum_u\left[\left(d_u-2\sum_v d_v\right) + d_u^2 -d_u^2 -\sum_v(d_u+1-3d_v) \right]
        \\ &=\frac{1}{n}\sum_u\left[\sum_v (d_v - d_u)\right] = 0
\end{align*}
because $\sum_u\sum_vd_v = \sum_u d_u^2$.
Going further,
\[
\begin{split}
    b_3 = \frac{1}{n}\sum_u\Bigg[\left(\sum_v\frac{1}{2} (3d_v^2-3d_v+2) + 3\sum_w d_wd_{uw}\right) + 2d_u\left(d_u-2\sum_v d_v\right) + d_u^3 \quad \\-d_u\left(\left(d_u-2\sum_v d_v\right)+d_u^2\right) - \sum_v(d_u+1-3d_v)d_u \quad\\- \frac{1}{2}\left(3+d_u-d_u^2-10d_v-6d_ud_v+16d_v^2\right)\Bigg]
\end{split}
\]
Using the fact that $w\in N^2(u)$ if and only if $u\in N^2(w)$, and then that $d_{uw}=d_{wu}$, one can show that 
\begin{align*}
    b_3 &= -\frac{1}{2n}\sum_u\left[d_u + 7\sum_v(d_u-d_v)^2\right] \le -\frac{1}{2}.
\end{align*}
This means that the desired inequality holds at the level of terms of order at most $\lambda^3$. 
Conceptually, this is perhaps enough to believe that subject to \emph{some} upper bound on $\lambda$ we have the desired result. We work a little harder to get a specific bound in terms of $\Delta$, although we do not compute a suitable absolute constant.

Quantitatively, it suffices to obtain the bound $\xi \le \lambda^3/2$ and to justify the expansion of $1/(1-y_u)$ by showing that $|y_u|\le 1/2$. 
These bounds follow from a crude by analysis of the expressions for $x_u$ and $y_u$. 
Since each of the degree variables $d_u$, $d_v$, $d_w$, and $d_{uw}$ that we can see in the sum are between $1$ and $\Delta$, it is easy to see that $x_u$ and $y_u$ are both bounded above and below by polynomials in $\Delta\lambda$ of degree at most four whose coefficients are absolute constants.
Crucially, these polynomials have no constant term. 
Then there is a positive absolute constant $c$ such that when $\lambda\le c/\Delta$ we have $|y_u|\le 1/2$. 
Turning to $\xi$, we see that it is bounded above by $P(\Delta\lambda)$ where $P(z)$ is a polynomial of degree $3\times 5 + 4 = 19$ in $z$ that is divisible by $z^4$. 
Then for $\lambda\le1/\Delta$ there is a positive absolute constant $c'$ such that $\xi \le c'(\Delta\lambda)^4$, and for another positive $c''$, when $\lambda \le c''/\Delta^4$ we have $\xi\le \lambda^3/2$ as required. 
\end{proof}

\section{Bounds on variance fraction}\label{sec:VAR}

In this section we formulate a higher-order version of the local occupancy methods whose roots are in~\cite{She95,Alo96,DJPR17a,DJPR17b}. 
The basic idea is to consider marginals or other probabilities as variables in an optimization problem whose objective is a quantity such as the occupancy or variance fraction. 
Using local information about the graph class being studied, one can find constraints satisfied by these variables which hold for any graph in the class. 
Typically, the objective and constraints are linear, and then one can consider the linear relaxation of this graph-defined optimization problem. 
Here, we extend this method in a new direction by considering the variance fraction as an objective, which is quadratic in such variables. 
One of the principal difficulties of this method is that finding constraints is not straightforward. 
We have fairly general methods for finding constraints but we do not have a classification that provides certainty, or even evidence, that one has exhausted the power of one's choice of variables for the optimization. 
Even worse, we do not have a framework for understanding when the linear (or now perhaps quadratic) relaxation of the optimization problem has a gap. 

\begin{proof}[Proof of Theorem~\ref{thm:Vbounds}]
    For convenience, fix $\lambda>0$ and let $r = \lambda/(1+n\lambda)$, $s=\lambda/(1+\lambda)$, and $\alpha=E_G(\lambda)$.
    Writing $p_u$ for the probability that $u$ is occupied and $p_{uv}$ for the probability that both $u$ and $v$ are occupied when an independent set is drawn from $\mu_{G,\lambda}$, we have 
    \begin{align*}
        V_G(\lambda) &= \frac{1}{n}\sum_{u}\left(p_u + \sum_{v\ne u}p_{uv} - p_u\sum_{v}p_v\right)
        \\ &= \alpha - n\alpha^2 + \frac{1}{n}\sum_u\sum_{v\ne u}p_{uv}.
    \end{align*}
    That is, the variance fraction is quadratic in the occupancy fraction and linear in a variable corresponding to some average of the pair marginals.

    For all vertices $u$ we have $0\le p_u\le s$.
    For all pairs $uv$ such that $u\ne v$ we have $0 \le p_{uv}\le sp_v\le s^2$ because $p_{uv}$ is a probability and 
    \[ p_{uv} = \Pr(u\in I \mid v\in I)\Pr(v\in I) \le sp_v. \]
    This holds because $s$ is an upper bound on $\Pr(u\in I \mid v\in I)$. 
    Note that the above bound is tight if and only if $u$ is isolated in $G-N[v]$.
    It also holds by Theorem~\ref{thm:Ebounds} that $r \le \alpha\le s$.
    
    Then we have the following upper bound. 
    \begin{align*}
        V_G(\lambda) &\le \alpha - n\alpha^2 + \frac{1}{n}\sum_u\sum_{v\ne u}sp_v
        \\&\le \alpha - n\alpha^2 + \frac{1}{n}\sum_us(n\alpha -p_u)
        \\&\le \alpha - n\alpha^2 +s(n-1)\alpha.
    \end{align*}
    This is increasing as a function of $\alpha$ for
    \[ \alpha \le \frac{1+s(n-1)}{2n}, \]
    and hence when $s\le 1/(n+1)$ the maximum occurs at the boundary: $\alpha=s$.  
    That is, under the condition $\lambda\le 1/n$ we have $s\le 1/(n+1)$ and that our upper bound on the variance fraction is an increasing function of $\alpha$. 
    At $\alpha=s$ we get the required bound $V_G(\lambda)\le s(1-s)$. 
    
    Under similar conditions, the minimization is analogous. 
    We have the lower bound. 
    \begin{align*}
        V_G(\lambda) &\ge \alpha - n\alpha^2.
    \end{align*}
    This is increasing when $\alpha\le 1/(2n)$ and hence if $s\le 1/(2n)$ we minimize the variance by minimizing $\alpha$ and setting it to $r$. We get 
    \[ V_G(\lambda) \ge r - nr^2 = \frac{\lambda}{(1+n\lambda)^2}, \]
    which is the desired lower bound. 
\end{proof}

\subsection{Examples of large variance fraction}

Some algebra and calculus reveals that for a $5$-vertex path $P_5$, we have 
\[ V_{P_5}(\lambda) > \frac{\lambda}{(1+\lambda)^2}\]
for all $\lambda \ge 33$. 

The family of cycle graphs $C_n$ is relatively straightforward to study because
\[ Z_{C_n}(\lambda) = 2(-\lambda)^{n/2} T_n\left(1/\sqrt{-4\lambda}\right), \]
where $T_n$ is the standard Chebyshev polynomial of the first kind given by $T_n(\cos\theta)=\cos(n\theta)$. 
This means that 
\[ V_{C_n}(\lambda) = \frac{n \sec ^2\left(n \sec ^{-1}\left(2 \sqrt{-\lambda}\right)\right)}{16\lambda+4}-\frac{\lambda \tan \left(n \sec ^{-1}\left(2 \sqrt{-\lambda}\right)\right)}{(-4\lambda-1)^{3/2}}, \]
and this example shows that 
\[ V_G(\lambda) / \frac{\lambda}{(1+\lambda)^2}\] 
is unbounded as $\lambda\to\infty$ (one can check that after taking the limit $n\to\infty$ this quantity is asymptotic to $\sqrt{\lambda}/8$).

\section{Free volume and related inequalities}\label{sec:FV}

We start with a nontrivial example of both FV and VAR inequalities holding. This example is inspired by the hard-core model on the smallest connected graph that is not a disjoint union of cliques: $K_{1,2}$.
We have to take a disjoint union of three copies of $K_{1,2}$ so that there is a disjoint union of cliques with the same degree sequence. Then we consider $P(x) = Z_{K_{1,2}}(x)^3$, and $Q(x) = Z_{K_2}(x)^3Z_{K_3}(x)$.

\begin{lemma}
Let 
\begin{align*}
    P(x) &= 1 + 9x + 30x^2 + 45x^3 + 30x^4 + 9x^5 + x^6 = (1+3x+x^2)^3\\
    Q(x) &= 1 + 9x + 30x^2 + 44x^3 +24x^4 = (1+2x)^3(1+3x).
\end{align*}
Then $P \geFV Q$ and $P\geVAR Q$.
\end{lemma}
\begin{proof}
    Straightforward computation suffices. To show $P \geFV Q$ requires checking five simple inequalities, and to show $P\geVAR Q$ we observe that 
    \[ P(x)^2Q(x)^2(V_P(x)-V_Q(x)) = 3 x^3 (2 x+1)^4 (x^2+3 x+1)^4 (86 x^4+176 x^3+118 x^2+32x+3),\]
    which is clearly nonnegative for all $x\ge 0$.
\end{proof}

A simple modification of $Q$ yields an example where VAR holds but not FV.
\begin{lemma}
Let 
\begin{align*}
    P(x) &= 1 + 9x + 30x^2 + 45x^3 + 30x^4 + 9x^5 + x^6\\
    Q(x) &= 1 + 9x + 30x^2 + 44x^3 +24x^4 + 9x^5.
\end{align*}
Then $P \not \geFV Q$ but $P\geVAR Q$.
\end{lemma}
\begin{proof}
The disproof of FV is straightforward. Note that 
\[ \frac{9}{30} = \frac{3}{10} \not\ge \frac{3}{8} = \frac{9}{24}. \]
The proof of VAR follows from the observation that 
\[ P(x)^2Q(x)^2(V_P(x)-V_Q(x))\]
is equal to 
\[\begin{split}
513 x^{21}+8136 x^{20}+56373 x^{19}+223380 x^{18}+554664 x^{17}+887748x^{16} + 926244 x^{15} + {}\\ 739620 x^{14}+864699 x^{13}+1472832 x^{12}+1950525x^{11}+1748244 x^{10}+1074906 x^9+ {}\\ 460512 x^8+ 137304 x^7+27864 x^6+3651x^5+276 x^4+9 x^3,
\end{split}\]
which is a polynomial in $x$ with nonnegative coefficients, and thus nonnegative for all $x\ge0$. 
\end{proof}
\begin{lemma}
Let 
\begin{align*}
    P(x) &= 1 + 9x + 30x^2 + 45x^3 + 30x^4 + 9x^5 + x^6\\
    Q(x) &= 1 + 9x + 30x^2 + 44x^3 +24x^4 + 10x^5.
\end{align*}
Then $P \not \geCOEF Q$ but $P\geVAR Q$.
\end{lemma}
\begin{proof}
    Clearly, COEF fails for the $x^5$ term. For VAR we used a computer algebra system to check symbolically that $V_P(x)\ge V_Q(x)$ for all $x\ge 0$.
\end{proof}

\begin{lemma}
    Let $P,Q$ be any of the following pairs of polynomials:
\begin{align*}
    P(x) &= 1 + 4x + 2x^2 + 2x^3 \\
    Q(x) &= 1 + 2x + x^2 + x^3,\\[\baselineskip]
    P(x) &= 1 + 10x + 210x^2 + 21x^3 + 21x^4 + 21x^5\\
    Q(x) &= 1 + 10x + 10x^2 + x^3 + x^4 + x^5,
\shortintertext{and}
    P(x) &= 1 + 10x + x^2 + 20010x^3 + 2001x^4 + 2001x^5\\
    Q(x) &= 1 + 10x + x^2 + 10x^3 + x^4 + x^5.
\end{align*}
Then $P \geFV Q$ but $P\not\geVAR Q$. 
\end{lemma}
\begin{proof}
    The FV inequalities are easy to check by hand. 
    For the first pair we have $V_Q(1) = \frac{26}{25} > \frac{74}{81} = V_P(1)$.
    For the second pair we have $V_Q(1) = \frac{53}{48} > \frac{18619}{20164} = V_P(1)$. For the third pair we have $V_Q(1) = \frac{293}{192} > \frac{68604293}{192384192} = V_P(1)$.
\end{proof}

The reason we give these counterexamples is to show that some additional properties that would be true in the setting of the hard-core model on two graphs of the same order can be satisfied, namely that the first few coefficients of the polynomials are equal.

\subsection{Counterexamples to upper bounds in the degree sequence setting}

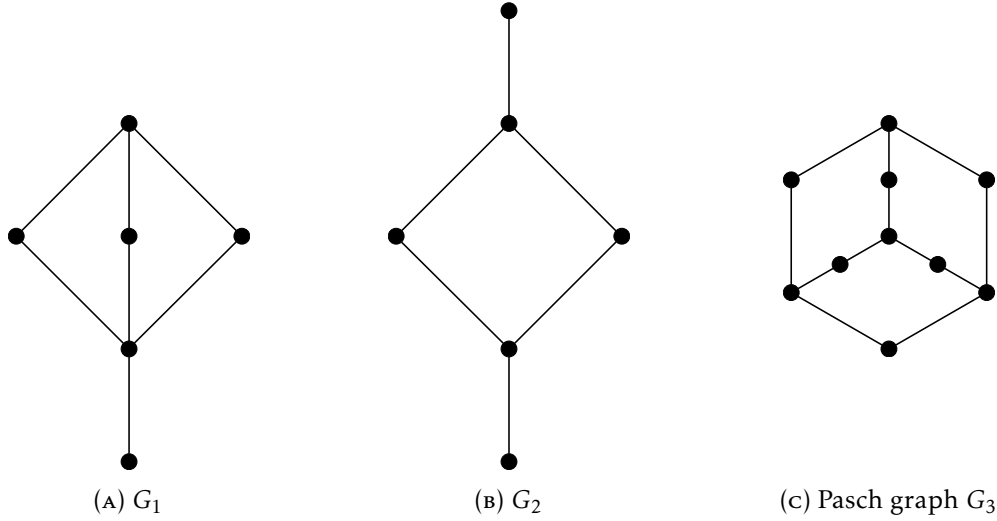
\begin{figure}[htb]
    \centering
    \subcaptionbox{$G_1$\label{fig:6,32}}[0.3\textwidth]{    \begin{tikzpicture} 
        \GraphInit[vstyle=Empty]
        \SetVertexSimple[MinSize=8pt,LineWidth = 0pt,LineColor = black,FillColor = black]
        \SetVertexNoLabel
        \Vertex[x=2cm,y=4cm]{v0}
        \Vertex[x=0cm,y=2cm]{v1}
        \Vertex[x=2cm,y=0cm]{v2}
        \Vertex[x=4cm,y=2cm]{v3}
        \Vertex[x=2cm,y=2cm]{v4}
        \Vertex[x=2cm,y=-2cm]{v5}
        \Edge[](v0)(v1)
        \Edge[](v1)(v2)
        \Edge[](v2)(v3)
        \Edge[](v3)(v0)
        \Edge[](v0)(v4)
        \Edge[](v4)(v2)
        \Edge[](v2)(v5)
    \end{tikzpicture}
    \begin{tikzpicture} 
        \GraphInit[vstyle=Empty]
        \SetVertexSimple[MinSize=8pt,LineWidth = 0pt,LineColor = black,FillColor = black]
        \SetVertexNoLabel
        \Vertex[x=2cm,y=4cm]{v0}
        \Vertex[x=0cm,y=2cm]{v1}
        \Vertex[x=2cm,y=0cm]{v2}
        \Vertex[x=4cm,y=2cm]{v3}
        \Vertex[x=-2cm,y=2cm]{v4}
        \Vertex[x=6cm,y=2cm]{v5}
        \Edge[](v0)(v1)
        \Edge[](v1)(v2)
        \Edge[](v2)(v3)
        \Edge[](v3)(v0)
        \Edge[](v1)(v4)
        \Edge[](v3)(v5)
    \end{tikzpicture}
    \begin{tikzpicture} 
        \GraphInit[vstyle=Empty]
        \SetVertexSimple[MinSize=8pt,LineWidth = 0pt,LineColor = black,FillColor = black]
        \SetVertexNoLabel
        \foreach\i in {0,1,...,5}
        {
            \pgfmathtruncatemacro{\a}{\i*60 + 30}
            \Vertex[a=\a,d=2cm]{v\i}
        }
        \Vertex[a=0,d=0]{vc}
        \foreach\i in {6,7,8}
        {
            \pgfmathtruncatemacro{\a}{\i*120 + 90}
            \Vertex[a=\a,d=1cm]{v\i}
        }
        \Edge[](v0)(v1)
        \Edge[](v1)(v2)
        \Edge[](v2)(v3)
        \Edge[](v3)(v4)
        \Edge[](v4)(v5)
        \Edge[](v5)(v0)
        \Edge[](vc)(v6)\Edge(v6)(v1)
        \Edge[](vc)(v7)\Edge(v7)(v3)
        \Edge[](vc)(v8)\Edge(v8)(v5)
    \end{tikzpicture}}
    \subcaptionbox{$G_2$\label{fig:Matchstick6,20}}[0.3\textwidth]{    \begin{tikzpicture} 
        \GraphInit[vstyle=Empty]
        \SetVertexSimple[MinSize=8pt,LineWidth = 0pt,LineColor = black,FillColor = black]
        \SetVertexNoLabel
        \Vertex[x=2cm,y=4cm]{v0}
        \Vertex[x=0cm,y=2cm]{v1}
        \Vertex[x=2cm,y=0cm]{v2}
        \Vertex[x=4cm,y=2cm]{v3}
        \Vertex[x=2cm,y=2cm]{v4}
        \Vertex[x=2cm,y=-2cm]{v5}
        \Edge[](v0)(v1)
        \Edge[](v1)(v2)
        \Edge[](v2)(v3)
        \Edge[](v3)(v0)
        \Edge[](v0)(v4)
        \Edge[](v4)(v2)
        \Edge[](v2)(v5)
    \end{tikzpicture}
    \begin{tikzpicture} 
        \GraphInit[vstyle=Empty]
        \SetVertexSimple[MinSize=8pt,LineWidth = 0pt,LineColor = black,FillColor = black]
        \SetVertexNoLabel
        \Vertex[x=2cm,y=4cm]{v0}
        \Vertex[x=0cm,y=2cm]{v1}
        \Vertex[x=2cm,y=0cm]{v2}
        \Vertex[x=4cm,y=2cm]{v3}
        \Vertex[x=-2cm,y=2cm]{v4}
        \Vertex[x=6cm,y=2cm]{v5}
        \Edge[](v0)(v1)
        \Edge[](v1)(v2)
        \Edge[](v2)(v3)
        \Edge[](v3)(v0)
        \Edge[](v1)(v4)
        \Edge[](v3)(v5)
    \end{tikzpicture}
    \begin{tikzpicture} 
        \GraphInit[vstyle=Empty]
        \SetVertexSimple[MinSize=8pt,LineWidth = 0pt,LineColor = black,FillColor = black]
        \SetVertexNoLabel
        \foreach\i in {0,1,...,5}
        {
            \pgfmathtruncatemacro{\a}{\i*60 + 30}
            \Vertex[a=\a,d=2cm]{v\i}
        }
        \Vertex[a=0,d=0]{vc}
        \foreach\i in {6,7,8}
        {
            \pgfmathtruncatemacro{\a}{\i*120 + 90}
            \Vertex[a=\a,d=1cm]{v\i}
        }
        \Edge[](v0)(v1)
        \Edge[](v1)(v2)
        \Edge[](v2)(v3)
        \Edge[](v3)(v4)
        \Edge[](v4)(v5)
        \Edge[](v5)(v0)
        \Edge[](vc)(v6)\Edge(v6)(v1)
        \Edge[](vc)(v7)\Edge(v7)(v3)
        \Edge[](vc)(v8)\Edge(v8)(v5)
    \end{tikzpicture}}
    \subcaptionbox{Pasch graph $G_3$\label{fig:Pasch}}[0.3\textwidth]{    \begin{tikzpicture} 
        \GraphInit[vstyle=Empty]
        \SetVertexSimple[MinSize=8pt,LineWidth = 0pt,LineColor = black,FillColor = black]
        \SetVertexNoLabel
        \Vertex[x=2cm,y=4cm]{v0}
        \Vertex[x=0cm,y=2cm]{v1}
        \Vertex[x=2cm,y=0cm]{v2}
        \Vertex[x=4cm,y=2cm]{v3}
        \Vertex[x=2cm,y=2cm]{v4}
        \Vertex[x=2cm,y=-2cm]{v5}
        \Edge[](v0)(v1)
        \Edge[](v1)(v2)
        \Edge[](v2)(v3)
        \Edge[](v3)(v0)
        \Edge[](v0)(v4)
        \Edge[](v4)(v2)
        \Edge[](v2)(v5)
    \end{tikzpicture}
    \begin{tikzpicture} 
        \GraphInit[vstyle=Empty]
        \SetVertexSimple[MinSize=8pt,LineWidth = 0pt,LineColor = black,FillColor = black]
        \SetVertexNoLabel
        \Vertex[x=2cm,y=4cm]{v0}
        \Vertex[x=0cm,y=2cm]{v1}
        \Vertex[x=2cm,y=0cm]{v2}
        \Vertex[x=4cm,y=2cm]{v3}
        \Vertex[x=-2cm,y=2cm]{v4}
        \Vertex[x=6cm,y=2cm]{v5}
        \Edge[](v0)(v1)
        \Edge[](v1)(v2)
        \Edge[](v2)(v3)
        \Edge[](v3)(v0)
        \Edge[](v1)(v4)
        \Edge[](v3)(v5)
    \end{tikzpicture}
    \begin{tikzpicture} 
        \GraphInit[vstyle=Empty]
        \SetVertexSimple[MinSize=8pt,LineWidth = 0pt,LineColor = black,FillColor = black]
        \SetVertexNoLabel
        \foreach\i in {0,1,...,5}
        {
            \pgfmathtruncatemacro{\a}{\i*60 + 30}
            \Vertex[a=\a,d=2cm]{v\i}
        }
        \Vertex[a=0,d=0]{vc}
        \foreach\i in {6,7,8}
        {
            \pgfmathtruncatemacro{\a}{\i*120 + 90}
            \Vertex[a=\a,d=1cm]{v\i}
        }
        \Edge[](v0)(v1)
        \Edge[](v1)(v2)
        \Edge[](v2)(v3)
        \Edge[](v3)(v4)
        \Edge[](v4)(v5)
        \Edge[](v5)(v0)
        \Edge[](vc)(v6)\Edge(v6)(v1)
        \Edge[](vc)(v7)\Edge(v7)(v3)
        \Edge[](vc)(v8)\Edge(v8)(v5)
    \end{tikzpicture}\vspace*{1.5cm}}
    \caption{Selected counterexamples to an OCC strengthening of the degree-sequence upper bound in~\eqref{eq:Fdegrees}.}
    \label{fig:OCCcx}
\end{figure}

Sah, Sawnhey, Stoner and Zhao already showed~\cite{SSSZ19} that a natural vertex-based upper bound on the free energy density is false. 
That is, taking a vertex-based viewpoint as in the true lower bound in~\eqref{eq:Fdegrees}, one might guess that for all $G$ and all $\lambda>0$ we have 
\[ F_G(\lambda) \le \frac{1}{n}\sum_{u\in V(G)} F_{K_{d_u,d_u}}(\lambda), \]
but this is false at $\lambda=1$ for a 4-vertex path. 
They established the edge-based upper bound in~\eqref{eq:Fdegrees}, but we point out here that the natural strengthening of \emph{that} bound to the occupancy fraction cannot hold. 
We demonstrate two counterexamples on 6 vertices and another, the so-called Pasch graph\footnote{\url{https://mathworld.wolfram.com/PaschGraph.html}}, on 10 vertices. The latter is notable for having no vertices of degree one. See Figure~\ref{fig:OCCcx}. The Pasch graph is the easiest of the given examples on which to verify our calculations because every edge joins a vertex of degree two to a vertex of degree three. 
Note that none of these counterexamples involve vertices of degree zero, so we can drop the sum over such vertices in the desired bound (see~\eqref{eq:Fdegrees}).
The edge-based upper bound on occupancy fraction \emph{does} hold in the case of $\Delta$-regular graphs~\cite{DJPR17b}, and the Pasch graph shows that this is extremely fragile as measured by e.g.\ $\max_{uv\in E(G)}|d_u-d_v|$ in the sense that as soon as this is nonzero the bound can not hold.

For the graph $G_1$ in Figure~\ref{fig:6,32}, we have 
\begin{align*}
E_{G_1}(\lambda) &= \frac{\lambda (4 \lambda ^3+12 \lambda ^2+16 \lambda +6)}{6(\lambda ^4+4 \lambda ^3+8 \lambda ^2+6 \lambda +1)},
\end{align*}
\begin{align*}
\frac{1}{6}\sum_{uv\in E(G_1)} &\frac{d_u+d_v}{d_ud_v}E_{K_{d_u,d_v}}(\lambda)= \\&\frac{1}{6} \left(\frac{\lambda  \left(3 (\lambda +1)^2+2 (\lambda+1)\right)}{2 \left((\lambda +1)^3+(\lambda +1)^2-1\right)}+\frac{\lambda\left(4 (\lambda +1)^3+1\right)}{4 \left((\lambda +1)^4+\lambda\right)}+\frac{3 \lambda  \left(4 (\lambda +1)^3+2 (\lambda +1)\right)}{8\left((\lambda +1)^4+(\lambda +1)^2-1\right)}\right),
\end{align*}
and the latter is smaller for e.g.\ $\lambda\ge 5$.

For the graph $G_2$ in Figure~\ref{fig:Matchstick6,20}, we have 
\begin{align*}
E_{G_2}(\lambda) &= \frac{\lambda (4 \lambda^3+12 \lambda^2+18 \lambda+6)}{6 (\lambda^4+4 \lambda^3+9 \lambda^2+6\lambda+1)},\\
\frac{1}{6}\sum_{uv\in E(G_2)} \frac{d_u+d_v}{d_ud_v}E_{K_{d_u,d_v}}(\lambda)&=\frac{1}{6} \left(\frac{2 \lambda \left(3 (\lambda+1)^2+1\right)}{3\left((\lambda+1)^3+\lambda\right)}+\frac{2 \lambda \left(3 (\lambda+1)^2+2 (\lambda+1)\right)}{3\left((\lambda+1)^3+(\lambda+1)^2-1\right)}\right),
\end{align*}
and the latter is smaller for e.g.\ $\lambda\ge 5$.

For the Pasch graph $G_3$ in Figure~\ref{fig:Pasch}, we have 
\begin{align*}
E_{G_3}(\lambda) &= \frac{\lambda \left(6 \lambda^5+30 \lambda^4+80 \lambda^3+126 \lambda^2+66 \lambda+10\right)}{10 \left(\lambda^6+6\lambda^5+20 \lambda^4+42 \lambda^3+33 \lambda^2+10 \lambda+1\right)},\\
\frac{1}{10}\sum_{uv\in E(G_3)} \frac{d_u+d_v}{d_ud_v}E_{K_{d_u,d_v}}(\lambda)&=\frac{\lambda  \left(3 (\lambda +1)^2+2 (\lambda +1)\right)}{5 \left((\lambda+1)^3+(\lambda +1)^2-1\right)},
\end{align*}
and the latter is smaller for e.g.\ $\lambda\ge 5$.

\subsection*{Acknowledgment}
We thank Will Perkins for comments on a draft of this work.

\addtocontents{toc}{\protect\setcounter{tocdepth}{-1}}

{
\small\hypersetup{urlcolor=MidnightBlue}
\bibliographystyle{habbrv}
\bibliography{main}
}

\addtocontents{toc}{\protect\setcounter{tocdepth}{1}}

\end{document}